\documentclass{article}

\usepackage{amssymb}
\usepackage{amsmath}
\usepackage{amsthm}

\newtheorem{theorem}{Theorem} \newtheorem{lemma}{Lemma}[section]
\newtheorem{propo}{Proposition}[section]
\newtheorem{defin}{Definition}

\def\Fo{F$\mbox{\o}$lner}
  \newcommand{\ep}{\varepsilon}
  
 \newcommand{\e}{\ep} 
\newcommand{\N}{\mathbb{N}}
\newcommand{\R}{\mathbb{R}} \newcommand{\Z}{\mathbb{Z}}

\newcommand{\G} {\mathcal{G}} 
 
\def\proof{\smallskip\noindent{\it Proof.} }

\newtheorem*{thma}{Theorem 1}
\newtheorem*{thmb}{Theorem 2}
\newtheorem*{thmc}{Theorem 3}

\newcommand{\urd} {U^r_d}

\newcommand{\omu} {\widehat{\mu}}
\newcommand{\ws} {\widehat{s}}\newcommand{\wt} {\widehat{t}}

\title{Hyperfinite actions on countable sets and probability measure spaces
\footnote{
\, Research sponsored by OTKA Grant No. 69062}}

\author{Mikl\'os Ab\'ert and G\'abor Elek}

\begin{document}

\title{Hyperfinite actions on countable sets and probability measure spaces}
\author{Mikl\'os Ab\'ert and G\'abor Elek \thanks{%
AMS Subject Classification: 43A07, 05C99 \, Research sponsored by OTKA Grant
No. 69062}}
\maketitle

\begin{abstract}
We introduce the notion of hyperfiniteness for permutation actions of countable groups  and give
a geometric and analytic characterization, similar to the known characterizations for amenable 
actions. We also answer a question of van Douwen on actions of the free group on two generators 
on countable sets.
\end{abstract}

\section{Introduction}

Let $\Gamma $ be a countable group acting on a countably infinite set $X$ by
permutations. An invariant mean on $X$ is a $\Gamma $-invariant, finitely
additive map $\mu $ from the set of subsets of $X$ to $[0,1]$ satisfying $%
\mu (X)=1$. Von Neumann \cite{neu} initiated the study of invariant means of
group actions.

We say that a group action of $\Gamma $ on $X$ is \emph{amenable} if $X$
admits a $\Gamma $-invariant mean. The group $\Gamma $ is amenable, if the
right action of $\Gamma $ on itself is amenable. Over the decades, amenability of groups has
become an important subject, with connections to various areas in
mathematics, like combinatorial group theory, probability theory, ergodic
theory and harmonic analysis.

All actions of amenable groups are amenable and for free actions, this
trivially goes the other way round as well, but in general, one has to
assume certain faithfulness conditions to make the notion meaningful. Even
when making the natural assumption that the action is transitive, the
general picture is that for most sets of conditions, one can construct
corresponding amenable actions of groups that are very far from being
amenable themselves.

In particular, van Douwen \cite{Dou} constructed a transitive amenable
action of the free group on two generators such that every nontrivial
element of the group fixes only finitely many points. We call this condition 
\emph{almost freeness}. Further examples of amenable actions of non-amenable
groups were given by Glasner and Monod \cite{GN} and by Moon \cite{Moon1},
\cite{Moon2}.

For probability measure preserving (p.m.p.) actions, the notion that mostly
corresponds to amenability is hyperfiniteness. Let $\Gamma $ act on a
probability measure space, preserving the measure. The action is called \emph{%
hyperfinite} if the measurable equivalence relation generated by the action
is up to measure zero an ascending union of finite measurable 
equivalence relations. As before,
all p.m.p. actions of amenable groups are hyperfinite, and for free actions,
there is equivalence, but in general, very large groups can act
hyperfinitely. In particular, Grigorchuk and Nekhrashevych \cite{Grig}
constructed an ergodic, faithful, hyperfinite p.m.p. action of a
non-amenable group. More generally, for a hyperfinite p.m.p. action of a
group $\Gamma $, the action of $\Gamma $ on almost all orbits is amenable.
For the other direction, Kaimanovich \cite{kai} presented a counterexample. 

The main goal of this paper is to introduce and analyze the notion of
hyperfiniteness for permutation actions of countable groups. If $\Gamma $
acts on a countably infinite set $X$, then the action always extends to the
Stone-Cech compactification $\beta X$. This connection establishes a
bijection between invariant means on $X$ and invariant measures on $\beta X$%
. In particular, an action is amenable if and only if the extended action
preserves a regular Borel-probability measure. This suggests the following
definition.

\begin{defin}
Let the countable group $\Gamma $ act on the set $X$ by permutations. We say
that the action is \emph{hyperfinite} if $\beta X$ admits a regular Borel
probability measure that is invariant under the extended action and for
which this action is hyperfinite. 
\end{defin}

In particular, every hyperfinite action is automatically amenable.

Our first result gives a combinatorial and a geometric characterization of
hyperfiniteness for actions of finitely generated groups. Let $G_{n}$ be a
sequence of graphs with an absolute bound on the degrees of vertices in $%
G_{n}$. We say that $(G_{n})$ is \emph{hyperfinite}, if for all $\varepsilon
>0$, there exists $Y_{n}\subseteq V(G_{n})$ and $K>0$ such that $\left\vert
Y_{n}\right\vert <\varepsilon \left\vert G_{n}\right\vert $ and every
connected component of the subgraph induced on $V(G_{n})\backslash Y_{n}$
has size at most $K$ ($n\geq 1$). This notion was introduced in \cite%
{Elekcost}.

\begin{thma}
Let $\Gamma $ be a group generated by the finite symmetric set $S$, acting
on the countably infinite set $X$ by permutations. Let $S_{\Gamma }$ denote the
Schreier graph of this action with respect to $S$. Then the following are
equivalent: \newline
1) The action is hyperfinite; \newline
2) There exists a hyperfinite \Fo-sequence in $S_{\Gamma }$; \newline
3) There exists an invariant mean $\mu $ on $X$, such that for all $%
\varepsilon >0$, there exists $Y\subseteq X$ with $\mu (Y)<\varepsilon $ and 
$K>0$ such that the connected components of the induced subgraph of $%
S_{\Gamma }$ on $X\backslash Y$ have size at most $K$.
\end{thma}

In \cite{Dou} van Douwen asked the following question. Let $H$ be a
countable infinite amenable group. Is there an almost free transitive action
of $F_{2}$, the free group of two generators, on $H$ such that every
invariant mean on $H$ is $F_{2}$-invariant? We will show that the answer is
negative, however, it is true if we change the almost freeness condition to
faithfulness. 

\begin{thmb}
\mbox{}
\begin{enumerate} 
\item   There exists no almost-free transitive action of $F_2$
on a finitely generated amenable group $H$ which preserves all 
$H$-invariant means.
\item  For any finitely generated amenable group $H$, there exists a faithful, 
transitive action of $F_2$ on $H$ which
preserves all the $H$-invariant means.
\end{enumerate}
\end{thmb}

Finally, we show the following. 

\begin{thmc}
There exists an ergodic, faithful p.m.p. profinite 
action of a non-amenable group that
is hyperfinite but topologically free. 
\end{thmc}

\noindent
Note that this answers a question of Grigorchuk, Nekrashevich and Sushchanskii
\cite{grineksus}. Note that Bergeron and Gaboriau \cite{GB} also constructed
an ergodic, faithful p.m.p profinite action which is not free, but
topologically free.

\section{The Stone-Cech compactification} \label{stone}
Let $X$ be a countably infinite set and $\beta X$ be its Stone-Cech
compactification. The elements of $\beta X$ are the ultrafilters on $X$
and the set $X$ itself is identified with the principal ultrafilters.
For a subset $A\subseteq X$, let $U_A$ be the set of ultrafilters
$\omega\in\beta X$ such that $A\in\omega$. Then $\{U_A\}_{A\subseteq X}$ forms
a base of the compact Hausdorff topology of $\beta X$. It is well-known that
the Banach-algebra of continuous functions 
$C(\beta X)$ can be identified with the
Banach-algebra $l^\infty(X)$. 
Let $\mu$ be a finitely additive measure on $X$. Then one can associate a
regular Borel measure $\omu$ on $\beta{X}$, by taking
$$\omu(U_A)=\mu(A)\,.$$
Indeed,  let $f\in l^\infty(X)$ be a bounded real function on $X$.
Then the continuous linear transformation 
$$T(f):=\int_{X} f \,d\mu$$ is well-defined. Thus, by the Riesz representation
theorem
$$T(f)=\int_{\beta X} f\,d\omu$$ for some regular Borel-measure $\omu$.
Since $T(\chi_{U_A})=\mu(A)$, the equality $\mu(A)=\omu(U_A)$ holds.
In fact, there is a one-to-one correspondance between the regular
Borel-measures and the finitely additive measures on $X$, since the integral
$\int_{X} f\,d\mu$ is completely defined by $\mu$.

\noindent
If $s:X\to X$ is a bijection, then it extends to a map
$\ws:\beta X\to\beta X$ by
$$\ws(\omega)=\bigcup_{A\in\omega} s(A)\,.$$
Since $\ws(U_A)=U_{s(A)}$, the map $\ws$ is a continuous bijection.
Thus if $\Gamma$ is a countable group acting on $X$, then we have an extended
action on $\beta X$. The following lemma is due to Bl\"umlinger \cite
[Lemma 1]{Blu}
\begin{lemma}
There is a one-to-one correspondance between the $\Gamma$-invariant means on $X$and the $\Gamma$-invariant regular Borel probability measures on $X$.
\end{lemma}
\proof
Observe that the set of $\Gamma$-invariant regular measures is the annihilator of
the set
$$\{f-\gamma(f)\,|\, f\in C(\beta X),\gamma\in \Gamma \}\,,$$
and the set of $\Gamma$-invariant means is the annihilator of the set
$$\{f-\gamma(f)\,|\, f\in l^\infty(X) \}\,.\quad\qed$$ 
Therefore an action of $\Gamma$ is amenable if
and only if the corresponding action on $\beta X$ has an invariant probability
measure.

\section{Geometrically hyperfinite actions}
Let $\Gamma$ be a finitely generated group acting on $X$ preserving the
mean $\mu$. Let $S$ be a finite, symmetric generating set for $\Gamma$ and
$S_\Gamma$ be the Schreier graph of the action. That is
\begin{itemize}
\item $V(S_\Gamma)=X$\,.
\item $(x,y)\in E(S_\Gamma)$ if $x\neq y$ and 
there exists $s\in S$ such that $s(x)=y$ 
\end{itemize}
Note that we do not draw loops in our Schreier-graphs.
Let $T$ be a subgraph of $S_\Gamma$. The edge measure of $T$ is defined as
$$\mu_E(T)=\frac{1}{2}\int_{X} \deg_T(x)\,d\mu(x)\,,$$
where $\deg_T(x)$ is the degree of $x$ in $T$.
We say that the action is {\it geometrically hyperfinite} if for any $\e>0$, there exists
$K_\e>0$ and a subgraph $G_\e\in S_\Gamma$ such that $V(G_\e)=X$ and
$$\mu_E(S_\Gamma\backslash G_\e)<\e\,$$
and all the components of $G_\e$ have size at most $K_\e$. It is easy to see
that geometrical hyperfiniteness does not depend on the choice of the
generating system. Note however, that the geometric hyperfiniteness and the
hyperfiniteness of an action do depend on the choice of the invariant measure.
It is possible that for some invariant mean $\mu_1$ the action is hyperfinite
and for another invariant mean $\mu_2$ the action is not hyperfinite, only
amenable.

\noindent
The hyperfiniteness of a family of finite graphs was introduced in
\cite{Elekcost}. Let $\G=\{G_n\}$ be a family of finite graphs with vertex
degree bound $d$. Then $\G$ is called hyperfinite if for any $\e>0$ there
exists $K_\e>0$ such that for any $n\geq 1$ one can delete $\e|V(G_n)|$ edges
from $G_n$ to obtain a graph of maximum component size at most $K_\e$.
\begin{propo} \label{hyperf}
Let $S_\Gamma$ be the Schreier graph of an action of the finitely generated 
group
$\Gamma$ on $X$. Then the following two statements are equivalent.
\begin{enumerate}
\item $S_\Gamma$ contains a hyperfinite \Fo-sequence.
\item The action is geometrically hyperfinite with respect to some $\Gamma$-invariant mean 
$\mu$.
\end{enumerate}
\end{propo}
Recall that a \Fo-sequence of $S_\Gamma$ is sequence of induced subgraphs
$\{F_n\}^\infty_{n=1}$, where the isoperimetric constant $i(F_n)$ tends to
zero as $n$ tends to infinity. The isoperimetric constant of a finite subgraph
is the number of outgoing edges divided by the number of vertices.
\vskip 0.15in
\noindent
\proof
Suppose that $S_\Gamma$ has a hyperfinite \Fo-sequence $\{F_n\}^\infty_{n=1}$.
Let $G_n\subseteq F_n$ be induced subgraphs such that $\lim_{n\to\infty}
\frac {|V(G_n)|}{|V(F_n)|}=1\,.$ Then clearly $\{G_n\}$ is a hyperfinite \Fo-
sequence as well. Therefore, we can suppose that $\{F_n\}^\infty_{n=1}$
are vertex-disjoint subgraphs. Indeed, let $F_{n_1}$ be an element of the
\Fo- sequence such that
$$\frac{|V(F_{n_1}\backslash F_1)|}{|V(F_{n_1})|}>1-\frac{1}{10}\,.$$
Then let $F_{n_2}$ be an element such that
$$\frac{|V(F_{n_1}\backslash (F_1\cup F_{n_1})|}{|V(F_{n_2}|}>
 1-\frac{1}{100}\,.$$

Inductively, we can construct a hyperfinite \Fo-sequence consisting
of vertex-disjoint subgraphs. Now let $\omega$ be an ultrafilter on $\N$ and 
$\lim_\omega$ be
the corresponding ultralimit $\lim_\omega:l^\infty(\N)\to\R$.
Let 
$$\mu(A):=\lim_\omega \frac {|A\cap V(F_n)|}{|V(F_n)|}\,.$$
Then $\mu$ is an invariant mean and the action is geometrically hyperfinite 
with 
respect to $\mu$.

\noindent
Now let us suppose that $\mu$ is a $\Gamma$-invariant mean on $X$ and the action
is
geometrically hyperfinite with respect to $\mu$. 
Let $\{G_\e\}_{\e>0}$ be the subgraphs of $S_\Gamma$ as in the definition of
hyperfiniteness.
We need the following lemma.
\begin{lemma}\label{estimate}
Let $R\subseteq S_\Gamma$ be a subgraph of components of size at most $C$. 
Suppose
that the edge-density (number of edges divided by the number of vertices) in
each component is at least $\alpha$.
Then $\alpha \mu(V(R))\leq \mu_E(R)\,.$
\end{lemma}
\proof
We can write $R$ as a vertex-disjoint union $R=\cup^k_{i=1} R_i$, where
all the components of $R_i$ are isomorphic, having $l_i$ vertices and
$m_i$ edges. Let $S_i\subset V(R_i)$ be a set containing exactly one vertex
from each component.  We can even suppose that under the isomorphisms of the
components we always choose the same vertex.
Thus by the invariance of the mean, we have a partition
$$V(R_i)=\cup^{l_i}_{j=1} S^j_i\,,$$
where $S^1_i=S_i$, $\mu(S^j_i)=\frac{1}{l_i} V(R_i)$, and $S^j_i$ also has the
property that it contains one vertex from each component and 
under the isomorphisms of the
components, it always contains the same vertex.
Then
$$\mu_E(R_i)=\frac{1}{2}\sum^{l_i}_{j=1} d^j_i \mu(S^j_i)\,,$$
where $d^j_i$ is the degree in a component of $R_i$ at a vertex of $S^j_i$.

\noindent
This yields $\mu_E(R_i)=m_i\mu(V(R_i))/l_i$. Therefore
$$\mu_E(R)=\sum_{i=1}^k \mu_E(R_i)=\sum^k_{i=1}\frac{m_i}{l_i}\mu(V(R_i))
\geq \alpha \mu(V(R))\,. \quad\qed$$

\vskip 0.2in
\noindent
Now, pick a sequence $\e_1\geq \e_2\geq\dots$ such that \begin{equation} 
\label{e1}
\sum_{i=1}^\infty \sqrt{\e_i} <1\,. \end{equation}
Let $\delta>0$ be a real number and $G_\delta$ be a subgraph as above.
Let $S^\delta_i$ be the union of components of $G_\delta$ in which the edge
density of $S_\Gamma\backslash G_{\e_i}$ is at least $\sqrt{\e_i}$. By the
previous lemma, we have
$$\mu(V(S^\delta_i))\sqrt{\e_i}\leq \mu_E(S_\Gamma\backslash G_{\e_i})\,.$$
Hence $\mu(V(S^\delta_i))\leq \sqrt{\e_i}$. By (\ref{e1}),
for any $n\geq 1$, there exists $G'_\delta\subset G_\delta$, having the
following properties.
\begin{itemize}
\item $G'_\delta$ is a union of components of $G_\delta$.
\item $\mu(V(G'_\delta))>0$.
\item If $Z$ is a component of $G'_\delta$ then
the edge-density of $S_\Gamma\backslash G_{\e_i}$ inside $Z$ is less than 
$\sqrt{\e_i}$, for any $1\leq i \leq n$. That is, we can 
remove $\sqrt{\e_i} |V(Z)|$ edges from $Z$ to obtain
a graph of maximum component size $K_{\e_i}$.
\end{itemize}

\vskip 0.1in
\noindent
For  $\e>0$ let $W^\e_\delta\subset G_\delta$ be the union of components
$H$ such that the isoperimetric constant of $H$ is less than $\e$. By our
previous lemma, it is easy to see that for any fixed $\e>0$ we have
$$\lim_{\delta\to 0} \mu(W^\e_\delta)=1\,.$$
Therefore, for any $n\geq 1$ there exists $\delta_n$ and 
a component $H_n$
of $G_{\delta_n}$ such that 
\begin{itemize}
\item the isoperimetric constant of $H_n$ is less than $\frac{1}{n}$,
\item for any $1\leq i \leq n$ one can remove $\sqrt{\e_i}|V(H_n)|$ edges from
  $H_n$ to obtain a graph of maximum component size $K_{\e_i}$.
\end{itemize}
This implies that $\{H_n\}^\infty_{n=1}$ is a hyperfinite 
\Fo-sequence in $S_\Gamma$.
\qed

\section{Graphs and graphings}
Let $T$ be a countable graph of vertex degree bound $d$, such that
$V(T)=X$. Then there exists an action of a finitely generated group $\Gamma$ such
that $T$ is the (loopless) Schreier graph of the action.
Indeed, one can label the edges of $T$ with finitely many labels
$\{c_1,c_2,\dots, c_n\}$ in such a way that incident edges are labeled
differently. This way we obtain the Schreier graph of the $n$-fold free
product of $C_2$. If $\mu$ is a $\Gamma$-invariant mean on $X$ such that $T$ is the
Schreier graph of the action with respect
to a finite symmetric generating set $S\subset\Gamma$, then $\mu$ is an $H$-invariant mean for any other
action by a finitely generated group $H$ with the same Schreier graph.
Indeed, if $h\in H$ is a generator of $H$ and $A\subseteq X$, then $A$
can be written as a disjoint union
$$A=\bigcup_{s\in S} A_s\cup A_1\,,$$
where $h(x)=s(x)$ on $A_s$ and $h(x)=x$ on $A_1$. 
Therefore,
$$\mu(h(A))=\sum_{s\in S}\mu(s(A_s))+\mu(A_1)=\mu(A)\,.$$
Thus if $T$ is a graph on $X$ with bounded vertex degree, we can actually
talk about $T$-invariant means on $X$.
Let us consider a $\Gamma$-action on $X$ preserving the mean $\mu$ and the
extended $\Gamma$-action on $\beta X$ preserving the associated probability measure
$\omu$.
Let $\G$ be the graphing of this action on $\beta X$ (see \cite{Kech})
associated to a finite symmetric generating set $S$, that
is the Borel graph on $\beta X$ such that $(x,y)\in E(\G)$ if $x\neq y$ and
$s(x)=y$ for some generator $s$.
\begin{lemma}
The graphing $\G$ depends only on $T$, assuming $T$ is the Schreier graph for
the action.
\end{lemma}
\proof Let $H$ be another finitely generated group with generating system
$S'$, such that the Schreier graph of this action is $T$ as well.
It is enough to prove that for any $\omega\in\beta X$ and $s'\in S'$ either
$s'(\omega)=\omega$ or $s'(\omega)=s(\omega)$ for some $s\in S$. Observe that 
there exists $s\in S$ or $s=1$ such that
$$B=\{n\in X\,\mid\, s(n)=s'(n)\}\in\omega\,.$$
Let $A\in\omega$. Then $s(A\cap B)=s'(A\cap B)$. Hence $s(A)\in s'(\omega)$.
Thus, $s(\omega)=s'(\omega)$.\qed

\vskip 0.2in
\noindent
We denote the graphing associated to $T$ by $\G(T)$.
Note that if $S$ is a subgraph of $T$, such that
$V(S)=A\subseteq X$ then $\G(S)\subseteq\G(T)$ and all the
edges of $\G(S)$ are in between points of $U_A$.
Let us briefly recall the local statistics for graphings \cite{Elekfin}.
A rooted, finite graph of radius $r$ is a graph $H$, with a distinguished
vertex $x$ such that
$$\max_{y\in V(H)} d_H(x,y)=r\,.$$

Let $\urd$ denote the finite set of rooted finite graphs of radius $r$ with
vertex degree bound $d$ (up to rooted isomorphism).
If $T$ is a countable graph as above, let $A(T,H)\subseteq X$ be the set of 
points $x$  such that the $r$-neighborhood of $x$ on $T$ is rooted isomorphic
to $H$.
Similarly, let $A(\G(T),H)\in\beta X$ be the set of points $\omega\in\beta X$
such that the $r$-neighbourhood of $\omega\in\G(T)$ is rooted isomorphic to
$H$.
We need the labeled version of the above setup as well. Let $U^{r,n}_d$
denote the finite set of rooted finite graphs of radius $r$ with vertex
degree bound $d$, edge-labeled by the set $[n]=\{1,2,\dots,n\}$. Now
let us label the edges of $T$ by the set $[n]$ in such a way that incident
edges are labeled differently. Then this labeling induces a Schreier graph
of $T$ and thus a labeling of $\G(T)$ as well. Again, for $\tilde{H}\in 
U^{r,n}_d$ let $A(T,\tilde{H})\subseteq X$ be the set of 
points $n$  such that the $r$-neighborhood of $n$ on $T$ is rooted isomorphic
to $\tilde{H}$.
Similarly, let $A(\G(T),\tilde{H})\in\beta X$ be the set of points 
$\omega\in\beta X$
such that the $r$-neighbourhood of $\omega\in\G(T)$ is rooted isomorphic to
$\tilde{H}$. For $\tilde{H}\in U^{r,n}_d$, we denote by $[H]$ the underlying
unlabeled rooted graph in $U^r_d$.
The following proposition states that the local statistics of $T$ and $\G(T)$ 
coincide.
\begin{propo}
For any $r\geq 1$ and $H\in\urd$
$$\mu(A(T,H))=\omu(A(\G(T),H))\,.$$
\end{propo}
\proof
Let us partition $A(T,H)$ into finitely many parts
$$A(T,H)=\cup_{\tilde{H}\,, [\tilde{H}]=H} A(T,\tilde{H})$$
\begin{lemma}
The $r$-neighborhood of any $\omega$ in $U_{A(T,\tilde{H})}$ is rooted-labeled
isomorpic to $\tilde{H}$. \end{lemma}
\proof
Let $\gamma$ and $\delta$ be elements in the $n$-fold free product of $C_2$
with word-length at most $r$.
Suppose that $\gamma(x)=\delta(x)$ if $x\in A(T,\tilde{H})$. Then
$\gamma(A)=\delta(A)$ if $A\subset A(T,\tilde{H})$, thus
$\gamma(\omega)=\delta(\omega)$.
Now suppose that $\gamma(x)\neq \delta(x)$ if $x\in  A(T,\tilde{H})$.
Then $\gamma\delta^{-1}(A(T,\tilde{H}))\cap A(T,\tilde{H})=\emptyset\,.$ 
Therefore,
$\gamma(\omega)\neq \delta(\omega)$. Therefore, the rooted-labeled $r$-ball
around $\omega$ is isomorphic to $\tilde{H}$. \qed

\vskip 0.1in
\noindent
By the lemma,

$$\omu(A(\G(T),H))\geq\sum_{\tilde{H}\,, [\tilde{H}]=H} 
\omu(U_{A(T,\tilde{H})})=
\sum_{\tilde{H}\,, [\tilde{H}]=H}\mu(A(T,\tilde{H}))=\mu(A(T,H))\,.$$
Since 
$$\sum_{H\in\urd} \mu(A(T,H))= \sum_{H\in\urd} \omu(A(\G(T),H))=1\,$$
the proposition follows. \qed

\vskip 0.1in
The following lemma is an immediate consequence of our proposition.
\begin{lemma} \label{fontoska}
Let $W\subseteq T$ be a subgraph. Then for any $l>0$ the $\mu$-measure of
points that are
in some components of size $l$ is exactly the $\omu$-measure of points that
are
 in
some components of $\G(W)$ of size $l$.
\end{lemma}
Now we prove the main result of this section.
\begin{theorem}\label{tetel1}
Let $\Gamma$ be a finitely generated group acting on $X$ preserving the mean
$\mu$.
Then the Schreier graph $S_\Gamma$ is hyperfinite if and only if the extended
action
on $\beta X$ is hyperfinite. That is, the notions of 
geometric hyperfiniteness and hyperfiniteness
of $\Gamma$-actions are equivalent.
\end{theorem}
\proof
First let us suppose that $S_\Gamma$ is hyperfinite. Then $S_\Gamma=W\cup Z$, where all
the components of $W$ have size at most $K$ and $\mu_E(Z)\leq \e$.
Then $\G(S_\Gamma)=\G(W)\cup\G(Z)$, where all the components of $\G(W)$ have
size at most $K$ and $\omu_E(\G(Z))\leq \e$. 
Hence if $S_\Gamma$ is hyperfinite, then the extended action is hyperfinite.

\noindent
Now let us suppose that the extended action on $\beta X$ is hyperfinite.
Let $T$ be the Schreier graph of the $\Gamma$-action on $X$ and $\e>0$ be a real
number. A basic subgraph of $T$ is a graph $(A,B,t)$, where
\begin{itemize}
\item $A$ and $B$ are disjoint subsets of $X$,
\item $t$ is a bijection between $A$ and $B$ with graph contained in $T$,
\item $V(A,B,t)=X$,
\item $E(A,B,t)=\bigcup_{x\in A} (x,t(x))$.
\end{itemize}
Clearly, $T$ can be written as an edge-disjoint union
$$T=\bigcup^n_{i=1}(A_i,B_i,t_i)\,.$$
Then $\G(T)=\bigcup^n_{i=1}(U_{A_i},U_{B_i},\wt_i)\,$ where
$\wt_i$ is the extension of $t_i$ to $\beta X$.
Since $\G(T)$ is hyperfinite, there exists $W\subset\G(T)$, such that
$$W=\bigcup^n_{i=1}(L_i,M_i,\wt_i)\,,$$
where
\begin{itemize}
\item $L_i\subset U_{A_i}$ is Borel for any $1\leq i \leq n$,
\item $\sum^n_{i=1}\omu(U_{A_i}\backslash L_i)\leq\frac{\e}{2}$,
\item $\wt(L_i)=M_i\,.$
\item all the components of $W$ have size at most $K$\,.
\end{itemize}
\begin{lemma}
For any $\delta>0$, there exist sets $N_i\subseteq A_i$ such that
$$\sum^n_{i=1}\omu(L_i\triangle U_{N_i})<\delta\,.$$ \end{lemma}
\proof
Since $\omu$ is regular, there exist compact sets $\{C_i\}^n_{i=1}$
and open sets $\{U_i\}^n_{i=1}$ such that
\begin{itemize}
\item $C_i\subset L_i\subset U_i\subset U_{A_i}$,
\item $\sum^n_{i=1} \omu(U_i\backslash C_i)<\delta\,.$
\end{itemize}
Cover $C_i$ by finitely many base sets $U_{A_{i,j}}$ that are contained
in $U_i$. Since the finite union of base sets is still a base set the lemma
follows.\qed

\vskip 0.2in
The following lemma is straightforward.
\begin{lemma} Let $W$ be as above and
let $\{W_k\}^\infty_{k=1}$ be a sequence
of subgraphings of $\G(T)$ such that
$\lim_{k\to\infty} \omu_E(W_k\triangle W)=0\,.$
Then $\lim_{k\to\infty} \omu(Bad^{W_k}_K)=0$, where
$Bad^{W_k}_K$ is the set of points that are in
a component of $W_k$ of size larger than $K$.
\end{lemma}
By the two lemmas, we have a sequence of subgraphings
$$W'_k=\bigcup^n_{i=1}(U_{N_{i,k}}, U_{t_i(N_{i,k})},\wt_i)$$
such that
$$\lim_{k\to\infty} \omu(Bad^{W'_k}_K)=0$$
and
$$\lim_{k\to\infty} \sum^n_{i=1}\omu(L_i\triangle U_{N_i,k})=0\,.$$
Now let us consider the subgraphs $H_k\subset T$
$$H_k=\bigcup^n_{i=1}(N_{i,k}, t_i(N_{i,k}), t_i)\,.$$
Then by Lemma \ref{fontoska},
$$\lim_{k\to\infty} \mu(Bad^{H_k}_K)=0\,.$$
and
$$\limsup_{k\to\infty} \mu_E(T\backslash H_k)\leq \e\,.$$
This immediately shows that $T$ is geometrically hyperfinite.  \qed

\vskip 0.2in
We finish this section with a proposition that further underlines
the relation between hyperfinite p.m.p actions and hyperfinite
actions on countable sets.
\begin{propo} \label{underline}
Let $\Gamma$ be a finitely generated group acting hyperfinitely p.m.p on
a probability measure space. Then for almost all orbits, the corresponding
actions are hyperfinite as well.
\end{propo}
\proof We use the same idea as in the proof of Proposition \ref{hyperf}.
Let $\e_1> \e_2 >\dots$ be real numbers such that \begin{equation} \label{est}
\sum^\infty_{n=1}\sqrt{\e_n}<\frac{1}{2}\,. \end{equation}
We may suppose that the graphing $\G$ of our action on the probability
measure space $(Z,\mu)$ is the ascending union of subgraphings
$\G=\cup^\infty_{n=1} \G_n$ such that for all $n\geq 1$
\begin{itemize}
\item $\mu_E(\G\backslash \G_n)<\e_n$
\item all the components of $\G_n$ are finite.
\end{itemize}
We can also suppose that the action is ergodic, since in the
ergodic decomposition $\mu=\int \mu_t d\nu(t)$, almost all the $\mu_t$
are hyperfinite actions as well \cite{Kech}.
Let $X_n\subset Z$ be the set of  points $z$ such that
the isoperimetric constant of the component of $z$ is greater than
$\sqrt{\e_n}$. Then using the same estimate as in Lemma \ref{estimate}, one
can immediately see that \begin{equation} \label{elsobecs}
\mu(X_n)\leq \sqrt{\e_n}\,. \end{equation}
Now let $Y^k_n\subset Z$ be the set of points $y$ in $Z$ such that
the edge density of $\G\backslash\G_k$ in the component of $\G_n$
containing $y$ is greater than $\sqrt{\e_k}$. Then \begin{equation}
\label{masodikbecs}
\mu(Y^k_n)\leq \sqrt{\e_k}\,. \end{equation}
For $k\geq 1$ we define the set $A_{k}\subset Z$ the following
way. The point $z$ is in $A_{k}$ if
\begin{itemize}
\item The component of $\G_{k+1}$ containing $z$ has isoperimetric constant
not greater than $\sqrt{\e_{k+1}}$;
\item For any $1\leq i \leq k$ the edge-density of $\G\backslash\G_i$ in
the component of $\G_{k+1}$ containing $z$
is not greater than $\sqrt{\e_i}$.
\end{itemize}
By (\ref{est}), the measure of $A_{k}$ is not zero. Therefore by ergodicity,
almost every point of $Z$ has an orbit containing a point from each $A_{k}$.
Let $z\in Z$ be such a point. Let $F_k$ be the component of $G_{k+1}$ of
a point in the orbit of $z$. Then by the two conditions above 
 $\{F_k\}$ is a hyperfinite
\Fo-sequence. \qed

\section{On faithfulness}
If $\Gamma$ acts on the countably infinite set preserving the mean $\mu$, 
faithfulness
means that for any $1\neq\gamma\in \Gamma$ the fixed point set of $\gamma$ is
not the whole set. Glasner and Monod proved that for any countable group
$\Gamma$ the free product $\Gamma\star \Z$ can act on a countable set in an
amenable, transitive and faithful manner. One can see however, that in their
construction, if an element $\gamma$ is in the group $\Gamma$, 
then the fixed point 
set of $\gamma$ has mean one. We call a group action preserving a mean $\mu$
 {\it
  strongly faithful} if the fixed point set of any non-unit element has
$\mu$-measure less than one. 
\begin{propo} If a countable group $\Gamma$ admits an amenable, strongly
  faithful action on countable set, then the group is sofic. \end{propo}
\noindent
 (see
  \cite{Elekhyp} for the definition of soficity).

\noindent
\proof
Recall that such an action is called {\it essentially free}, if the 
fixed point set of any non-unit element has $\mu$-measure zero. It is proved
in \cite[Corollary 4.2]{Elekhyp} that any countable group with an 
amenable essentially free
action is sofic. Hence, the only thing remaining is to show the following
lemma.
\begin{lemma} Let $\Gamma$ be a countable 
group acting amenably and strongly faithfully
  on a countably infinite set $X$, preserving the mean $\mu$. Then $\Gamma$ admits an
  amenable, essentially-free action on a countably infinite set.
\end{lemma}
\proof
Let $K=\bigcup_{n=1}^\infty X^n$. Let us consider the product action of
$\Gamma$ on $X^n$. Define the mean
$\mu_2$ the following way. If $A\subset X\times X$ let
$$\mu_2(A)=\int_X \mu(\pi_1((A\cap (X,z) ))d\mu(z)\,,$$
where $\pi_1$ is the projection to the first coordinate.
Clearly, $\mu_2$ is preserved by the $\Gamma$-action. Inductively, we can
construct invariant means $\{\mu_n\}^\infty_{n=1}$ on
 the sets $\{X^n\}^\infty_{n=1}$. Now, let $\omega$ be a non-principal
ultrafilter on $\N$. Let us define a mean on $K$ the following way.
$$\nu(B)=\lim_\omega \mu_n(B\cap X^n)\,.$$
Then $\mu$ is a $\Gamma$-invariant mean on $K$.
Let $1\neq\gamma\in\Gamma$ and $F$ be the fixed point set of $\gamma$ in $X$.
The fixed point set of $\gamma$ in $X^n$ is exactly $F^n$, and obviously,
$$\mu_n(F^n)=(\mu(F))^n\,.$$
Hence, $\nu(\cup^{\infty}_{n=1}F^n) =0$ and the lemma follows. \qed

\section{On a problem of van Douwen}
In \cite{Dou} van Douwen asked the following question [Question 1.4]: 
If $H$ is
any countable infinite amenable group, then is there an almost free transitive
action of $F_2$ (the free group of two generators) on $H$ such that
every invariant mean on $H$ is $F_2$-invariant ?

\begin{theorem}\label{tetel2} 
\mbox{}
\begin{enumerate} 
\item   There exists no almost-free transitive action of $F_2$
on a finitely generated amenable group $H$ which preserves all 
$H$-invariant means.
\item  For any finitely generated amenable group $H$, there exists a faithful, 
transitive action of $F_2$ on $H$ which
preserves all the $H$-invariant means.
\end{enumerate}
\end{theorem}
\proof Let $Cay(H,S)$ be the Cayley-graph of the finitely generated group $H$
with respect to a symmetric generating system $S$. Suppose that $F_2$
acts almost freely on $H$. We separate two cases for the action.

\noindent
{\bf Case 1} There exists a \Fo-sequence $F_1, F_2,\dots $ in
$Cay(H,S)$ for which the following holds.
\begin{itemize}
\item $\{sF_n\cup tF_n \cup s^{-1} F_n \cup t^{-1} F_n \cup F_n\}^\infty_{n=1}$
are disjoint subsets, where $s,t$ are generators of $F_2$.
\item There exists $\e>0$ such that for any $n\geq 1$,
$\frac{|(sF_n\cup tF_n \cup s^{-1} F_n \cup t^{-1} F_n)
\backslash F_n|}{|F_n|}>\e$\,.
\end{itemize}

We define the $H$-invariant mean $\mu$ by
$$\mu(A):=\lim_\omega \frac{|A\cap F_n|}{|F_n|}\,,$$
where $\omega$ is a nonprincipal ultrafilter on $\N$ and $\lim_\omega$ is the
corresponding ultralimit.
We claim that $\mu$ is not preserved by the $F_2$-action.
Observe that for any $n\geq 1$ at least one of the following four
inequalities hold:
$|s F_n\backslash F_n|\geq \frac{\e}{4}|F_n|, 
|t F_n\backslash F_n|\geq \frac{\e}{4}|F_n|,
|s^{-1} F_n\backslash F_n|\geq \frac{\e}{4}|F_n|,
|t^{-1} F_n\backslash F_n|\geq \frac{\e}{4}|F_n|\,.$
Hence we can assume that for the set $A$ defined by
$$A =\{n\,\mid\, |s F_n\backslash F_n|\geq \frac{\e}{4}|F_n|\}\,,
$$ $A\in\omega$. Therefore
$$\mu(s\cup^\infty_{n=1} F_n)<1\quad\mbox{and} \quad 
\mu(\cup^\infty_{n=1} F_n)=1\,.$$
Therefore $\mu$ is not preserved by the $F_2$-action.

\noindent
{\bf Case 2}
If \Fo-sequences described in {\bf Case 1} do not exist, then any
\Fo-sequence is almost invariant under the $F_2$-action, that is
\begin{equation} \label{doueq1}
\lim_{n\to\infty} \frac{|(sF_n\cup tF_n \cup s^{-1} F_n \cup t^{-1} F_n)
\backslash F_n|}
{|F_n|}=0\,.
\end{equation}

Now let us fix a \Fo-sequence $\{G_n\}^\infty_{n=1}$ in $Cay(H,S)$.
By \cite[Proposition 4.1]{ESZAME}, $\{G_n\}^\infty_{n=1}$ is a hyperfinite
sequence. That is, for any $\e>0$ there exists $K_\e>0$ such that
one can remove $\e|V(G_n)|$ vertices and the incident edges in such a
way that in the resulting graph $G'_n$, all components have size
at most $K_\e$. By the counting argument applied in the proof of
Proposition \ref{hyperf}, it is easy to see that one can even suppose
that all the remaining components have isoperimetric constant at
most $\sqrt{\e}$ in $G_n$.
Now let us consider the following graph sequence $\{T_n\}^\infty_{n=1}$
 edge-labeled by the
set $\{s,t,s^{-1}, t^{-1}\}$
\begin{itemize}
\item
$V(T_n)=V(G_n)$
\item
$(x,y)$ is a directed edge labeled by $s$ (resp. by $t$, $s^{-1}$ or $t^{-1}$)
if $s(x)=y$ (resp. $t(x)=y$, $s^{-1}(x)=y$ or $t^{-1}(x)=y$)\,.
\end{itemize}
By almost-freeness and (\ref{doueq1}) it is clear that $\{T_n\}^\infty_{n=1}$
is a sofic approximation of $F_2$ (see \cite{ESZAME} for definition).

\begin{lemma}
$\{T_n\}^\infty_{n=1}$ is a hyperfinite graph sequence.
\end{lemma}
\proof
By (\ref{doueq1}) there exists a function $f:\R\to\R$ such
that $\lim_{x\to 0} f(x)=0$ satisfying
$$|(sL \cup tL \cup s^{-1}L \cup t^{-1}L)\backslash L |< f(\delta)|L|$$
for any finite set $L\subset H$ with isoperimetric constant less than
$\delta$.

\noindent
Let $\{G'_n\}^\infty_{n=1}$ be the subgraphs of $\{G_n\}^\infty_{n=1}$
obtained by removing $\e|V(G_n)|$ vertices and the incident edges
such that all the components of $G_n'$ have size at most $K_\e$ and
$G_n$-isoperimetric constant at most $\sqrt{\e}$.
The number of edges of $T_n$ that are in between the components of $G'_n$ is
less than $4 f(\sqrt{\e}) |V(T_n)|$. Hence by removing
$4 f(\sqrt{\e}) |V(T_n)|$ edges from $T_n$ and all the edges that are
incident to a vertex in $V(T_n)\backslash V(G_n')$ we can obtain a
graph with maximum component size at most $K_\e$. 
Therefore
$\{T_n\}^\infty_{n=1}$ is hyperfinite. \qed

\noindent
Since by \cite[Proposition 4.1]{ESZAME}, $F_2$ has no hyperfinite
sofic approximation, we obtain a contradiction.
Therefore $F_2$ has no almost-free action on $H$ that preserves all the
$H$-invariant means.

\vskip 0.2in
\noindent
Now we construct a faithful and transitive $F_2$-action on $H$ that
preserves all the $H$-invariant means.
First, fix a subset $\{i_n\}^\infty_{n=-\infty}\subset \Z$ such
that $i_n< i_{n+1}$ for all $n\in \Z$. Then fix a function
$f:\Z\to\{1, -1\}\,.$ The action $\alpha$ of $F_2$ on $\Z$ is defined the 
following way. Let
$s$ and $t$ be the generators of $F_2$.
\begin{itemize}
\item If $n$ is odd and $i_n<j< i_{n+1}$, then let $s(j)=j$.
\item If $n$ is even and $f(i_n)=1$, then if
$i_n\leq j < i_{n+1}$, let $s(j)=j+1$. If $j=i_{n+1}$, let
$s(j)=i_n$.
\item If $n$ is even and $f(i_n)=-1$, then if $i_n<j\leq i_{n+1}$, let
$s(j)=j-1$. If $j=i_n$, let $s(j)=i_{n+1}$.
\item If $n$ is even and $i_n<j< i_{n+1}$, then let $t(j)=j$.
\item If $n$ is odd  and $f(i_n)=1$, then if
$i_n\leq j < i_{n+1}$, let $t(j)=j+1$. If $j=i_{n+1}$, let
$t(j)=i_n$.
\item If $n$ is odd and $f(i_n)=-1$, then if $i_n<j\leq i_{n+1}$, let
$t(j)=j-1$. If $j=i_n$, let $t(j)=i_{n+1}$.
\end{itemize}

Note that the orbits of the $F_2$-action $\alpha$ generated by $s$ 
(we call these orbits
$s$-orbits) resp. by $t$ are finite cycles.
Clearly, one can define $\{i_n\}^\infty_{n=-\infty}$ and $f:\Z\to\{1,-1\}$
in such a way that for any $1\neq \gamma\in F_2$ there exists $n\in\Z$
such that $\gamma(n)\neq n$.

\noindent
Now let $\phi:\Z\to H$ be a bijection and $K>0$ such that
$d(\phi(n),\phi(n+1))\leq K$ for any $n\in \Z$, where $d(x,y)$ is the
shortest path distance in the Cayley graph of $H$.
Such bijection always exists from $\Z$ to an infinite connected bounded
degree graph $G$ if $G$ has one or two ends \cite{Seward}.
On the other hand, the Cayley graph of an amenable group always has one
or two ends \cite{MV}.

\noindent
The $F_2$-action on $H$ is given by
$$\gamma(x):=\phi(\alpha(\gamma)\phi^{-1}(x))\,.$$
Let $\mu$ be an $H$-invariant mean. We need to prove that $\mu$ is invariant
under the $F_2$-action above.

\begin{lemma} \label{cycle}
Let $n\geq 1$ and let
\begin{itemize}
\item $\Omega^s_n:=\{x\in H\,\mid\, d(s(x),x) \geq Kn\}$
\item $\Omega^{s^{-1}}_n:=\{x\in \Z\,\mid\, d(s^{-1}(x),x)\geq Kn \}$
\item $\Omega^t_n:=\{x\in \Z\,\mid\, d(t(x),x) \geq Kn \}$
\item $\Omega^{t^{-1}}_n:=\{x\in \Z\,\mid\, d(t^{-1}(x),x)\geq Kn\}$
\end{itemize}
Then the $\mu$-measure of any of these sets is less than $\frac{1}{n}\,.$
\end{lemma}
\proof Clearly, each $s$-orbit of size at least $n+1$
contains at most one element of $\Omega^s_n$.
The other $s$-orbits are disjoint from $\Omega^s_n$.
We need to show that the union of $s$-orbits of the $F_2$-action on $H$ of size
 at 
least $n+1$ has
measure at least $n \mu( \Omega^s_n)\,.$
Let $C$ be such an $s$-orbit of size $t$ and let $C_p$ be the unique
vertex such that
$$d(s^i(C_p),s^{i+1}(c_p))\leq K$$
for $i\leq t-2$.
Consider the set 
$$\bigcup_{C, |C| \geq n+1}\,\,C_p=L\,.$$
It suffices to prove that 

\begin{equation} \label{egye3} \mu(L)=\mu(s(L))=\dots=\mu(s^n(L))
\end{equation}
Define $h_p\in H$ by $s(C_p)=h_pC_p\,.$ Then $h_p$ is in
the $K$-ball around the unit element in the Cayley-graph of $H$.
Let $L=\bigcup_{h\in B_K(1)} L^h$, where $C_p\in L^h$ if $s(C_p)=hC_p$.
Then
$$ \mu(s(L))=\sum_{h\in B_k(1)} \mu(s(L^h))=
\sum_{h\in B_k(1)} \mu(L^h)=\mu(L)\,.$$
Similarly, $\mu(s^i(L))=\mu(L)$ if
$i\leq t-1$, therefore (\ref{egye3}) holds. \qed

\vskip 0.2in
\noindent
Now we finish the proof of the second part of our theorem.
Let $A\subseteq X$. Then
$$A=\bigcup_{h\in H} A_h\quad\mbox{where}\quad A_h=\{x\in A\,\mid\,s(x)=hx\}$$
Obviously, $\mu(A_h)=\mu(s(A_h))\,.$
Also, by our previous lemma, 
$$\mu(\bigcup_{h, h\notin B_{Kn}(1)}s( A_h))\leq\frac{1}{n}\,.$$
Therefore $\mu(s(A))\leq\mu(A) +1/n\,.$ for any $n\geq 1$. Hence
$\mu(s(A))\leq \mu(A)\,.$
However, the same way we can see that $\mu(s^{-1}(A))\leq \mu(A)$ as well.
That is, $\mu(s(A))=\mu(A)\,.$ Similarly, $\mu(t(A))=\mu(A)\,.$ \qed

\section{A topologically free, hyperfinite action of a nonamenable group}
Answering a question of Grigorchuk, Nekhrashevich and Sushschanskii 
\cite{grineksus}
Gaboriau and Bergeron \cite{GB} constructed a profinite, faithful, ergodic 
action of a
nonamenable group that is not essentially free, but topologically free. Note
that topological freeness of an action means that the set of points that are
not in the fixed point set of any nonunit element of the group is comeager.

On the
other hand, Grigorchuk and Nekrashevich constructed a profinite, ergodic
action of a nonamenable group that is faithful and hyperfinite. Their
construction is very far from being topologically free, in fact the set of
points that are not in the fixed point set of any nonunit element is meager.
However, we prove that the two results can be combined.
\begin{theorem} \label{tetel3}
There exists a finitely generated non-amenable group with an ergodic,
faithful, profinite action that is hyperfinite and topologically free.
\end{theorem}
\proof  An {\it amoeba} is a finite connected
 graph $G$ (with loops) having edge-labels
$A,B,C,D$  satisfying the following properties:
\begin{itemize}
\item $G$ is the union of simple cycles $\{C_i\}^{n_G}_{i=1}$. Some of the
  cycles might be loops. We call these cycles the basic cycles.
\item Any two of the basic cycles intersect each other in at most one vertex.
\item Let us consider the graph $T$ where the vertex set of $T$ is the set of
  basic cycles and two vertices are connected if and only if the corresponding
  cycles have non-empty intersection. Then $T$ is a tree.
\item For each vertex $x\in V(G)$ and for any label $A,B,C,D$  there exists
exactly one edge (maybe a loop) incident to $x$ having that label.
Hence the degree of any vertex is $4$. Note that the contribution of a loop
in the degree of a vertex is $1$.
\item Each loop is labeled by $C$ or $D$. For any vertex $x$ there are $0$ or 
$2$ loops incident to $x$.
\end{itemize}
Clearly, any amoeba is a planar graph. The minimal amoeba $M$ has $2$
vertices. The edge set of $M$ consists of a cycle of length two labeled by
$A$ and $B$ respectively and two loops for each of the two vertices labeled
by $C$ and $D$.

\noindent
Let $G$ be an amoeba. A {\it doubling} of $G$ is a two-fold topological
covering $\phi:H\to G$ by an amoeba $H$ constructed the following way.
Let $\{C_i\}^{n_G}_{i=1}$ be the set of basic cycles of $G$. One way to
construct $H$ is to pick a basic cycle $C_i$ which is not a loop and consider
a two-fold covering $\psi:C'_i\to C_i$. Obviously, $\psi$ extends to a
two-fold covering $\phi:H\to G$ uniquely, and $H$ is an amoeba as well.

\noindent
The second way to construct $H$ is to pick a vertex $x$ incident to two
loops $C_i$ and $C_j$. Then we cover their union by a cycle of length
two. Again, this covering extends to a two-fold covering in a unique way.
\subsection{The cycle-elimination tower}
Let start with a minimal amoeba $G_1$. A {\it cycle-elimination tower}
is a sequence of doublings
$$G_1\stackrel{\phi_1}{\leftarrow} G_2 \stackrel{\phi_2}{\leftarrow} 
G_3\stackrel{\phi_3}{\leftarrow}\dots$$
such that there exists a sequence of vertices $\{p_n\in
V(G_n)\}^\infty_{n=1}$,
 $\phi_n(p_{n+1})=p_n$ with the following property.
For any $k\geq 1$, there exists an integer $n_k$ such that the 
$k$-neighborhood of $p_{n_k}$ in the graph $G_{n_k}$ is a tree. It is
easy to see that by succesively eliminating cycles, such a tower can be
constructed.

\noindent
Let $\Gamma$ be the group $C_2\star C_2\star C_2\star C_2$ with free
generators $A,B,C,D$ of order $2$. Note that $\Gamma$ acts on the vertices
of an amoeba. Indeed, the generator $A$ maps the vertex $x$ to the unique
vertex $y$ such that the edge $(x,y)$ is labeled by $A$. If $x=y$, that is
the edge is a loop, then $A$ fixes $x$. Since the covering maps commute with
the $\Gamma$-actions, one can extend the $\Gamma$-action to the inverse limit 
space
$X_\Gamma=\lim_{\leftarrow} V(G_n)$. Recall that there is a natural probability
measure $\mu$ on $X_\Gamma$ induced by the normalized counting measures on the
vertex sets $V(G_n)$. 
The $\Gamma$-action on $\mu$ preserves the measure $\mu$ and
in fact this is the only Borel probability measure preserved by the action. 
The ergodicity of the $\Gamma$-action follows from the fact that $\Gamma$ acts
transitively on each vertex set $V(G_n)$.
\subsection{Topological freeness}
In this subsection we show that the action of $\Gamma$ on $X_\Gamma$ is
topologially free. Let us introduce some notation. If $m>n$, let
$\phi^m_n$ be the covering map from $G_m$ to $G_n$. Also,
let $\Phi_n:X_\Gamma\to G_n$ be the natural covering map from the inverse limit
space. We need to prove that if $1\neq \gamma\in \Gamma$, then
the fixed point set of $\gamma$ has empty interior.

Let $q\in V(G_n)$. Then $\Phi_n^{-1}(q)$ is an basic open set in $X_\Gamma$.
It is enough to prove that there exists $z\in \Phi_n^{-1}(q)$ such
that $\gamma(z)\neq z$. Let $d=dist_{G_n}(q,p_n)$, where $dist$ is the
shortest path distance and $\{p_n\}^\infty_{n=1}$ is the sequence of vertices
as above.
By the properties of graph coverings, for any
element $x$ in $(\phi^m_n)^{-1}(p_n)$ there exists $r\in (\phi^m_n)^{-1}(q)$
such that $dist_{G_m}(x,r)=d$. 

Now let $w(\gamma)$ be the wordlength of $\gamma$ and consider the
vertex $x=p_{n_{d+|w(\gamma)|}}\,.$
Clearly, if $r\in G_{n_{d+|w(\gamma)|}}$ and $dist(r,x)=d$ then
$\gamma(r)\neq r$. Therefore $\Phi^{-1}_n(q)$ contains a point $z\in X_\Gamma$
that is not fixed by $\gamma$.

\subsection{Hyperfiniteness}
Now we finish the proof of Theorem \ref{tetel3} by showing that the action
of $\Gamma$ on $X_\Gamma$ is hyperfinite.
Fix $\e>0$. Let us recall \cite{BSS} that planar graphs with bounded
vertex degree form a hyperfinite family. Hence there exists $K>0$
such that for each $G_n$ one can remove
$\frac{\e}{10}|V(G_n)|$ edges in such a way that in the remaining
graph $G'_n$ the maximal component size is at most $K$.

\noindent
Since $X_\Gamma$ is the inverse limit of $\{V(G_n)\}^\infty_{n=1}$ for 
any $p\in X_\Gamma$ there exists $m(p)\in\N$ such that if $l\geq m(p)$
then the $K+1$-neighborhood of $p$ in its $\Gamma$-orbit graph and the
$K+1$-neighborhood of $\Phi_l(p)$ in $G_l$ are isomorphic.

\noindent
Hence we have $m>0$ such that the Haar-measure of the set $A$ of the points
in $X_\Gamma$ for which the $K+1$-neighborhood of $x$ is
not isomorphic to the $K+1$-neighborhood of $\Phi_m(x)\in G_m$ is less than
$\frac{\e}{10}$. 

Let $\G$ denote the graphing of the $\Gamma$-action on $X_\Gamma$.
We remove the edges from $\G$ that are incident
to a point in $X$. Also, we remove the edges that are inverse images
of an edge removed from $E(G_n)$. Then the edge-measure of the edges
removed from $\G$ is less than $\e$ and in the remaining graphing
all the components have size at most $K$. \qed

\end{document}